\documentclass{amsart}
\usepackage{graphicx}
\usepackage{url}
\usepackage{subfigure}
\usepackage{amssymb,amsthm, amsmath}
\usepackage[foot]{amsaddr}

\def\N{\mathbb N}

\def\R{\mathbb R}

\def\GSE{\mathrm{GSP}}
\def\T{\mathrm{T}}
\def\P{\mathrm{\Pi}}
\def\PD{\mathrm{PD}}
\def\ZD{\mathrm{ZD}}
\def\ZT{\mathrm{ZT}}

\newtheorem{thm}{Theorem}

\newtheorem{lem}{Lemma}
\newtheorem{rem}{Remark}

\begin{document}
\title{An alternative proof for an aperiodic monotile}
\author{Shigeki Akiyama}
\address{Institute of Mathematics, University of Tsukuba}
\email{akiyama@math.tsukuba.ac.jp}
\author{Yoshiaki Araki}
\address{Japan Tessellation Design Association}
\email{yoshiaki.araki@tessellation.jp}
\date{}
\maketitle

\begin{abstract}
We give a simple alternative proof that 
the monotile introduced by \cite{SMKGS:23_1} is aperiodic.
\end{abstract}

\section{Introduction}
\label{Intro}
Smith Hat tiles 
found in \cite{SMKGS:23_1} form a one-parameter family Tile $(b:1)$ of aperiodic monotiles: each member tiles a plane but only non-periodically 
except when $b=0,1,\infty$.  
In this paper, we pick out one special case $b=\sqrt{3}$: the tile called 
Smith Turtle or simply Turtle, 
and give a simple alternative proof of this fact.
\begin{figure}[h]
\begin{center}
\includegraphics[width=6cm]{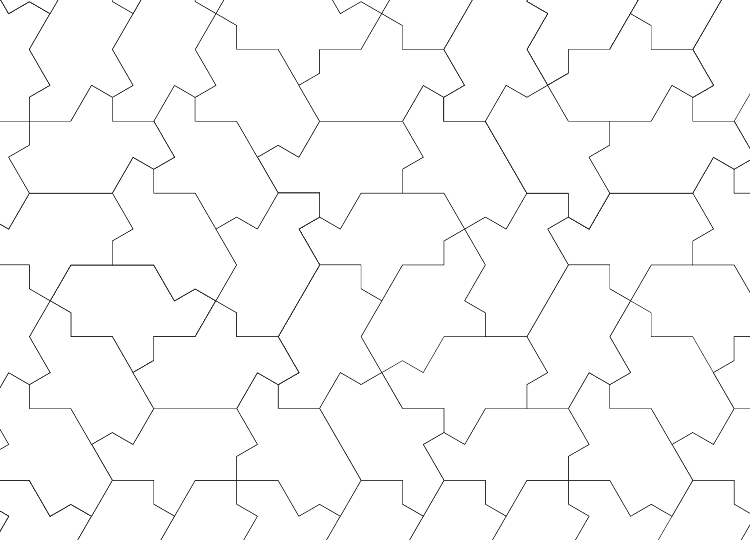}
\end{center}
\caption{Tiling by Smith Turtle}
\end{figure}

The proof in \cite{SMKGS:23_1} that it tiles a plane depends on a ``combinatorial'' substitution rule, which differs each time of its application and converges to a geometric substitution rule whose attractors give a limit tiling having fractal boundaries.
By this combinatorial nature, the proof becomes a little involved. 

In \S \ref{GoldenHex}, we give a concrete ``Golden Hex substitution",
whose tiles are essentially regular triangles and parallelograms,
but we additionally use approximate 
`linear' patches to fill gaps in the construction.
We call these linear patches ``Golden Sturmian Patches", 
see Figure \ref{GSPDef}.
Golden Sturmian Patches encode sturmian words of slope $(5-\sqrt{5})/10$.

We briefly recall
sturmian words and central words in \S \ref{Sturm}.
Central words are special palindromes that appear in sturmian words. 
They form a kind of building block of sturmian words and 
play an essential role in describing combinatorial properties
of sturmian words.
The consistency of the new substitution rule 
is shown by an induction, and the proof 
heavily rely on the property of central words. 
This construction of patch-tiles
proves that tilings by Smith Turtle do exist.

In \S \ref{GoldenAmmann}, we show that all tilings generated by Turtle are
non-periodic, using a special linear marking; ``Golden Ammann bar". We call it GAB in short. 
Such markings were 
originally introduced by R.~Ammann to construct 
a finite set of tiles having aperiodicity \cite{AGS:92,Gruenbaum-Shephard:87,Akiyama:12}. 
The choice $b=\sqrt{3}$ is essential to introduce
these Ammann bars to Tile $(b:1)$. 
Our proof requires only elementary properties of this GAB that 
are easily checked, 
and it does not require the meta-tiles or 
substitution structures found in \cite{SMKGS:23_1}.
%Since a tiling by Tile $(b:1)$ are uniquely transformed to the one by
%Tile $(\sqrt{3},1)$ and vice versa when $b\neq 0,1,\infty$, 

Tilings produced by Tile $(b:1)$ and the ones by  Tile $(c:1)$ are combinatorially equivalent if $b,c\not \in\{0,1,\infty\}$. 
Combining the above discussion,
we obtain a simple independent proof of the
aperiodicity of Smith Hat tile.
%The markings corresponds to the tiling generated by Golden hex substitution is non periodic.
%Indeed, the bars consists 
%a lacunary hexagonal grid whose gaps of adjacent lines
%give $1-2$ sequences, which forms sturmian sequences 
%of golden mean slope.
%We also discuss cut and project scheme associated with the fractal tilings using the algebraic dual scheme.
%We see that the fractal tilings generated by Tile $(b:1)$ does not depend on 
%$b$ up to similitudes and 
%the division of internal space gives a tiling consists of fractal curves and segments. 
%Indeed, we rediscover that the segments forms a Golden Hex tiling, which gives a self-duality of this aperiodic monotile. 

\section{Sturmian words and central words}
\label{Sturm}

Let $\{0,1\}^*$ be the set of finite words over $\{0,1\}$; a monoid by
concatenation operation equipped with the identity:
the empty word $\lambda$.
The set of right infinite words over $\{0,1\}$ is denoted by 
$\{0,1\}^{\N}$.
Sturmian words are the elements in $\{0,1\}^{\N}$ that
emerge from
codings of irrational rotations. 
Here we recall definitions and properties
that we shall use in this paper. 
Let $\alpha\in (0,1)$ be an irrational number and take
$\rho\in [0,1)$.
A sturmian word of slope $\alpha$ is an infinite word over $\{0,1\}$
defined by
$$
\lfloor \alpha(n+1) +\rho \rfloor - \lfloor \alpha n +\rho \rfloor
$$
or
$$
\lceil \alpha(n+1) +\rho \rceil - \lceil \alpha n +\rho \rceil
$$
for $n=1,2,3,\dots$. 
%A bi-infinite 
%Sturmian words in $\{0,1\}^{\Z}$ is obtained similarly 
%by taking index $n\in \Z$.
From the continued fraction expansion of
$$
\alpha=\cfrac
{1}{a_1+\cfrac {1}{a_2+\cfrac{1}{\ddots}}}=[a_1,a_2,\dots],
$$
{\bf standard words} $s_n$ for $n\ge -1$ can be defined by the recurrence
$$
s_{-1}=1,\ s_{0}=0,\ s_{n+1}=\overbrace{s_n s_n\dots s_n}^{d_n \text{times}} 
s_{n-1}=s_n^{d_n}s_{n-1}
$$
where $d_1=a_1-1$ and $d_n=a_n\ (n=2,3,\dots)$, 
see \cite[Proposition 2.2.24]{Lothaire:02}. 
The sturmian word of slope $\alpha$ and $\rho=0$ is obtained as the limit of $s_n$ as $n\to \infty$.
Every subword of the sturmian word of slope $\alpha$
is a subword of $s_n$ for some $n$.
A standard word $s$ of length greater than $1$ has three palindrome subwords $p,q,r$ that
$$
s=pw=qr\quad (w\in \{01,10\})
$$
and a standard word without specifying a slope 
is characterized by this property.
The decomposition $s=qr$ is unique if $s$ is not a palindrome, see
\cite[Theorem 2.2.4]{Lothaire:02}.
The palindrome $p$ is called the {\bf central word}. 
Central words are often called bi-special words, and literally
play the central role in describing the combinatorial and dynamical
properties
of the sturmian word (c.f.\ \cite[Chapter 6]{Fogg:02}).
Roughly speaking, a 
sturmian word is constructed by gluing together
central words with 
the `paste' words $\{01,10\}$.
Let $\tau=(\sqrt{5}+1)/2$ be the golden ratio
and fix $$\alpha=(5-\sqrt{5})/10=1/(1+\tau^2)
=[3,1,1,1,\dots]$$ and apply these formulas.
Then $s_0=0, s_1=001$, and $s_{n+1}=s_{n}s_{n-1}\ (n\ge 1)$
and we have
$$
s_2=0010,\ s_3=0010001,\ s_4=00100010010,\ s_5=001000100100010001, \dots
$$
Since $s_n$ is not a palindrome for $n\ge 1$, we 
define corresponding palindromes by $p_n,q_n,r_n$. Thus
$$
p_0=\lambda, \ 
p_1=0,\ p_2=00,\ p_3=00100,\ p_4=001000100,\ p_5=0010001001000100, \dots,
$$
and we obtain $s_{2n-1}=p_{2n-1}01,\ s_{2n}=p_{2n}10$ for $n\ge 1$.
Consequently 
$$
s_{2n-1}=s_{2n-2}s_{2n-3}=p_{2n-2}10 p_{2n-3}01$$
and
$$
p_{2n-1}=p_{2n-2}10 p_{2n-3},\  
q_{2n-1}=p_{2n-2},\ r_{2n-1}=10p_{2n-3}01
$$
hold. Similarly, we see
$$
s_{2n}=s_{2n-1}s_{2n-2}=p_{2n-1}01 p_{2n-2}10
$$
and
$$
p_{2n}=p_{2n-1}01 p_{2n-2},\ 
q_{2n}=p_{2n-1},\ r_{2n}=01p_{2n-2}10.
$$
Summarizing these, we have

\begin{equation}
\label{Cent1}
p_{2n+1}=p_{2n}10 p_{2n-1}=p_{2n-1}01p_{2n-2}10p_{2n-1}=
p_{2n-1}01p_{2n}
\end{equation}
and
\begin{equation}
\label{Cent2}
p_{2n}=p_{2n-1}01p_{2n-2}=p_{2n-2}10p_{2n-3}01p_{2n-2}=
p_{2n-2}10p_{2n-1}.
\end{equation}
\noindent
These decomposition rules are used in the next section.

\section{Golden Hex substitution}
\label{GoldenHex}

\begin{figure}[h]
\begin{center}
\includegraphics[width=15cm]{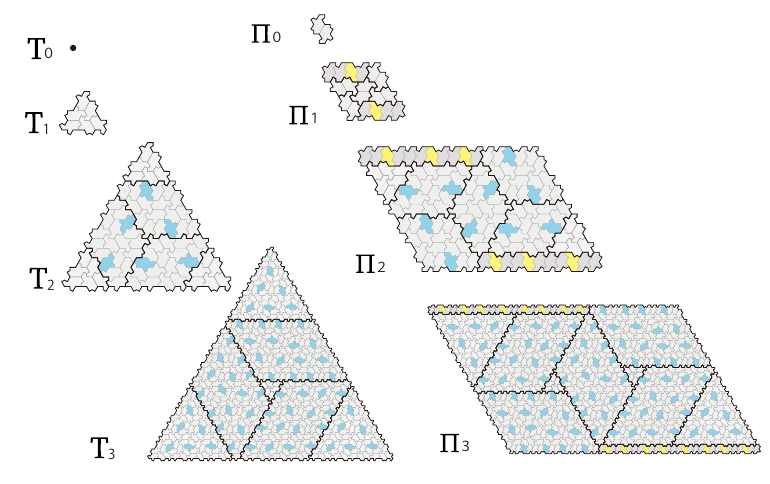}
\end{center}
\caption{Patch-tiles $\T_n, \P_n\ (n=0,1,2,3)$\label{GH2}}
\end{figure}

\begin{figure}[h]
\subfigure[$\T$]{%
\includegraphics[width=6cm]{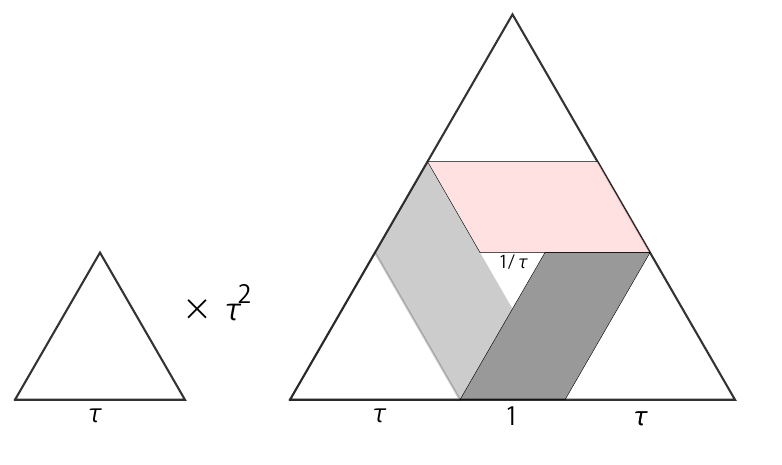}
}%
\subfigure[$\P$]{%
\includegraphics[width=8cm]{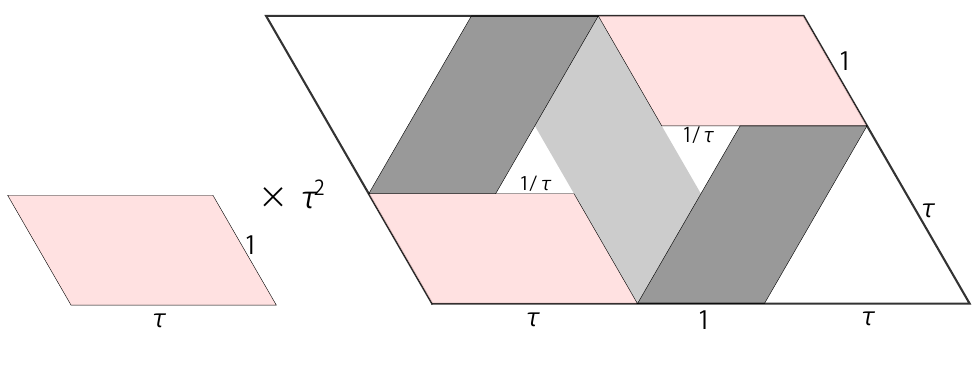}
}%
\caption{Golden Hex Substitution\label{GHex}}
\end{figure}

In this paper, we assume that a tile is a set homeomorphic to a closed ball.
A patch is a finite collection of tiles so that distinct tiles have disjoint interiors. A tiling is the covering of $\R^2$ by tiles where distinct tiles
have disjoint interiors, using 
only finitely many different tiles up to rigid motion.
Given a tiling, if we can 
partition it into a finite set of patches up to rigid motion so that
each of which is reconsidered as a new tile, then
we call them {\bf patch-tiles}, a similar idea is found 
in \cite{Frank-Sadun:14}.
We shall define two growing sequences $(\T_n,\P_n)_{n=0,1,\dots}$
of patch-tiles generated by Turtle depicted in Figure \ref{GH2}
\noindent
whose limit substitution is given by Figure~\ref{GHex}. 
The patch-tile $\T_n$ is invariant by $2\pi/3$-rotation 
and $\Pi_n$ is invariant by $\pi$-rotation.
Flipped tiles are 
colored blue or yellow.

\begin{figure}
\includegraphics[width=13cm]{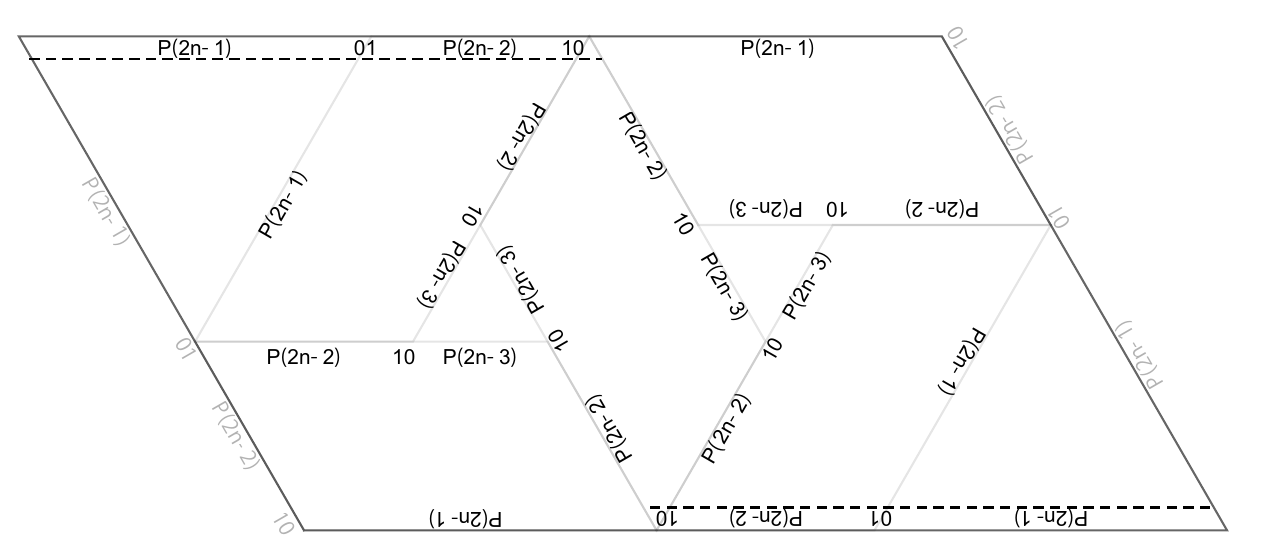}
\caption{$\P_n= 5 \P_{n-1} + 2 \T_{n-1} + 2 \T_{n-2} + 2 s_{2n}$\label{P}}
\end{figure}

We use an abusive terminology ``Golden Hex substitution", to refer to both this limit substitution in Figure \ref{GHex} and
the sequence of patch-tiles that approximates this substitution rule. Indeed, 
patch-tiles $\T_n$ and $\P_n$ are essentially regular
triangles and parallelograms, and $T_0$ is a single point. 
In the top left and bottom right parts of $\P_n\ (n=1,2,3)$ 
in Figure~\ref{GH2}, 
gray-yellow `linear' patches fill gaps to form the next level. 
This is made precise as broken lined parts in Figure~\ref{P}. 

\begin{figure}[h]
\begin{center}
\includegraphics[width=6cm]{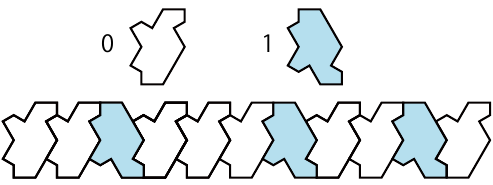}
\end{center}
\caption{Geometric realization of $00100010010$. This is a GSP.
\label{GSPDef}}
\end{figure}

For any word of $\{0,1\}^*$, we associate a geometric realization as 
a linear patch as in Figure~\ref{GSPDef}. 
Here the majority 
letter 0 corresponds to the orientation of the
Turtle whose head points 
upwards.
The geometric realization of the identity: empty word
$\lambda$ is set to be a single point.
In this paper, a {\bf Golden Sturmian Patch}, in short {\bf GSP}, 
is a geometric realization of a subword of 
the sturmian word of slope $(5-\sqrt{5})/10$. GSPs are the above
linear patches that fill the gaps to the next level, 
and give us an explicit construction of the sequence of patch-tiles. 

\begin{figure}[h]
\subfigure[$\T_n$]{%
\includegraphics[width=5cm]{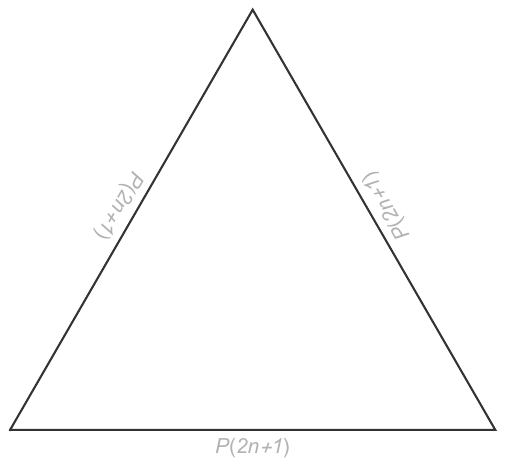}
}%
\subfigure[$\P_n$]{%
\includegraphics[width=7cm]{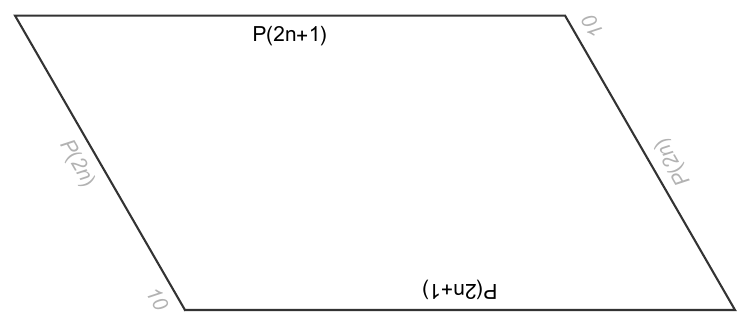}
}%
\caption{$\T_n$, $\P_n$ and surrounding GSPs \label{TP0}}
\end{figure}

\begin{figure}[h]
\includegraphics[width=9cm]{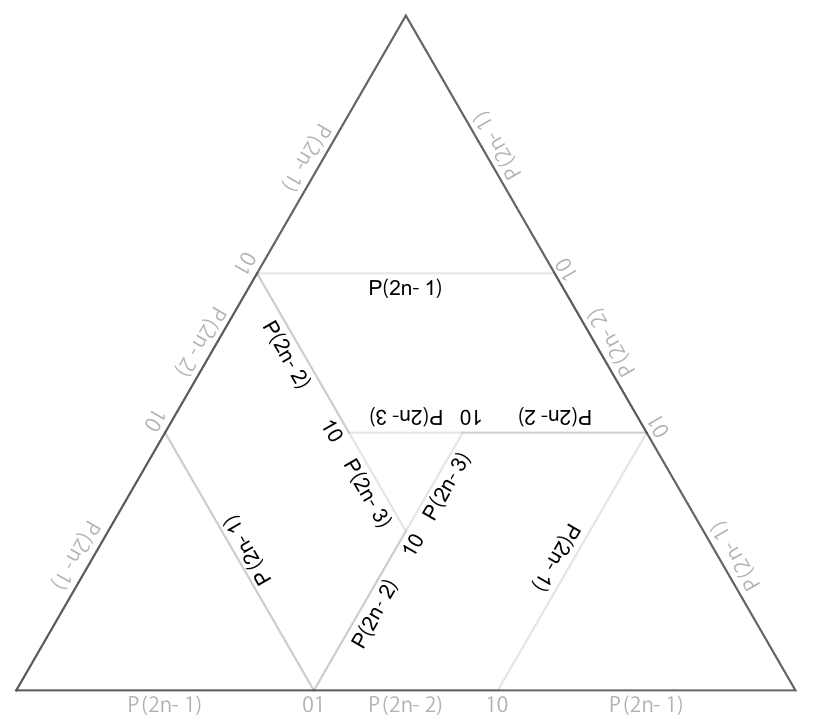}
\caption{$\T_n = 3 \T_{n-1} + 3 \P_{n-1} + \T_{n-2}$
\label{T}}
\end{figure}

For brevity, we denote the GSP of $p_n$ by $P(n)$, and write $01,10$ to express the GSP of $01,10$. Since they appear along with $P(n)$, 
there is no room for confusion.
The rotated $P(n)$ and $01,10$ are expressed by the corresponding rotated words by the same angle.
Then the inductive construction 
on $\T_n, \P_n$ is described in Figures~\ref{TP0},
\ref{T} and \ref{P}. 

Here $\T_n$ is subdivided into three $\T_{n-1}$'s, three
$\P_{n-1}$'s and one $\T_{n-2}$, and $\P_n$ is subdivided into five $\P_{n-1}$'s, two $\T_{n-1}$'s, two $\T_{n-2}$'s and two GSPs of $s_{2n}$ up to rotation.
The gray GSPs on the boundary mean it can receive the corresponding GSPs in the indicated direction. Figure \ref{TPlevel3} shows how the
gray GSPs work.

The role of GSPs in the substitution rule 
delayed our comprehension of the Golden Hex substitution. 
See the \S \ref{Appe}
for our previous understanding. 
Here is a key property in the proof of Theorem \ref{GHSTiling}.

\begin{figure}
\includegraphics[width=8cm]{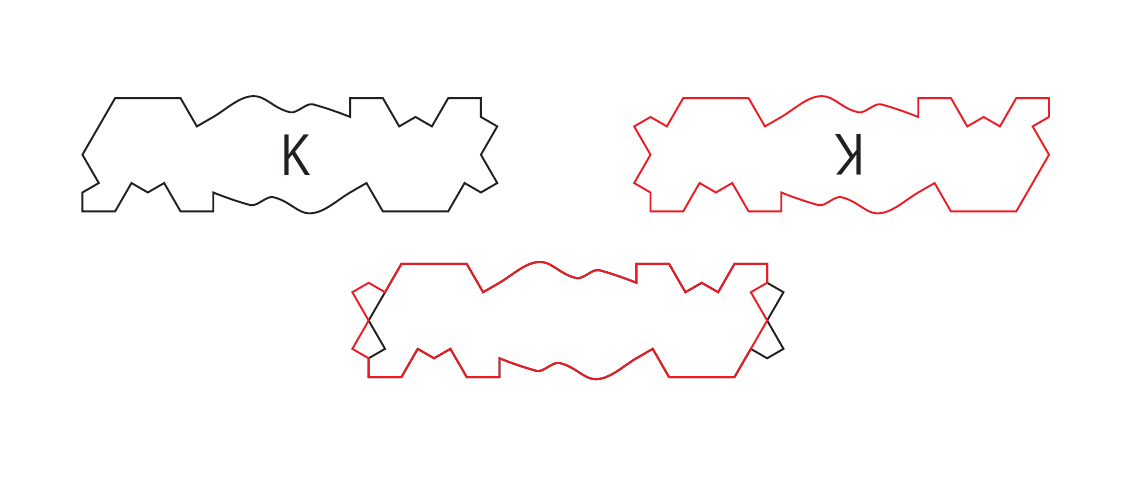}
\caption{A palindrome $K$ and its $\pi$-rotation are only different by 4 kites in both ends.\label{Pal4}}
\end{figure}

\begin{lem}
\label{Pal}
For any palindrome $K\in \{0,1\}^*$, the geometric realization of $K$ and its $\pi$-rotated image $\rotatebox[origin=c]{180}{K}$ differ only at the two ends by four small kites as in Figure~\ref{Pal4}.
In particular, $P(n)$ and $\rotatebox[origin=c]{180}{P(n)}$ share the same upper and lower boundaries.
\end{lem}

\proof
It is proved by a plain induction on the number of tiles, since in a palindrome of length $n\ge 2$, there is a palindrome subword 
of length $n-2$ in the middle. 
\qed

\begin{thm}
\label{GHSTiling}
The patch-tile sequence $(\T_n,\P_n)_{n=0,1,2,\dots}$ is well defined. 
\end{thm}
\noindent
Thus there are patches containing 
arbitrary large balls 
and the tiling by Turtle does exist (c.f.\ \cite[Section 3.8]{Gruenbaum-Shephard:87}).

\begin{proof}
We can check that $\T_n,\P_n$ with $n=0,1,2$
satisfies the condition of Figure~\ref{TP0}
as in Figure~\ref{Ind1}.\footnote{$\P_0$ and its surrounding in Figure~\ref{Ind1} is not used later in the proof.}
%, and probably
%this configuration around $\P_0$ does not appear in the tiling.}

\begin{figure}[h]
\includegraphics[width=8cm]{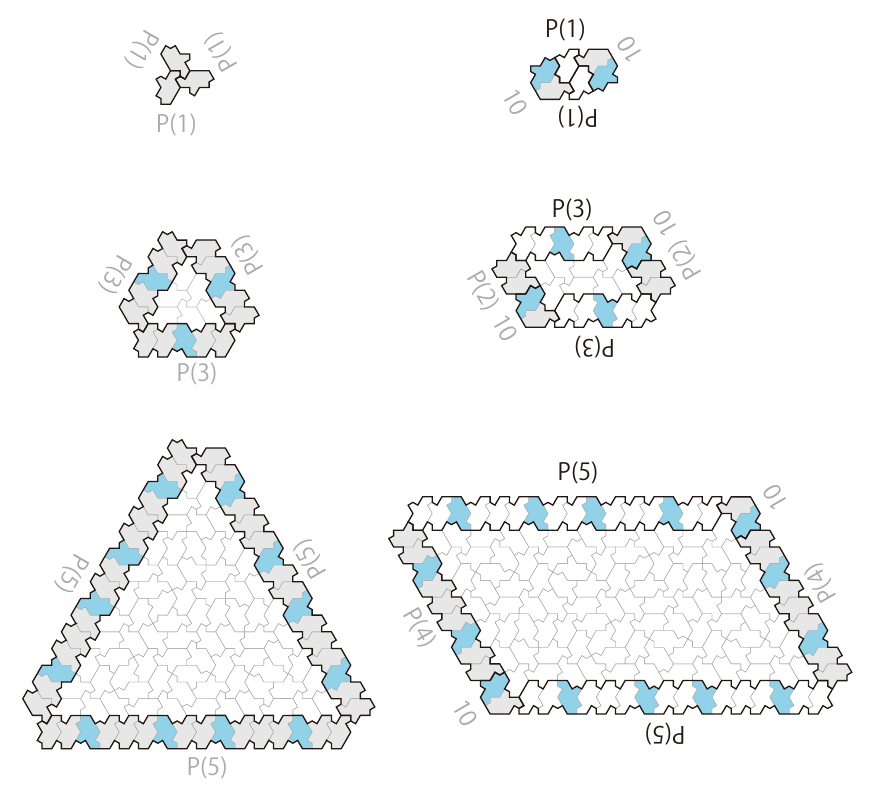}
\caption{The condition of $\T_i,\ \P_i\ (n=0,1,2)$ 
in Figure~\ref{TP0}\label{Ind1}}
\end{figure}

\begin{figure}
\includegraphics[width=15cm]{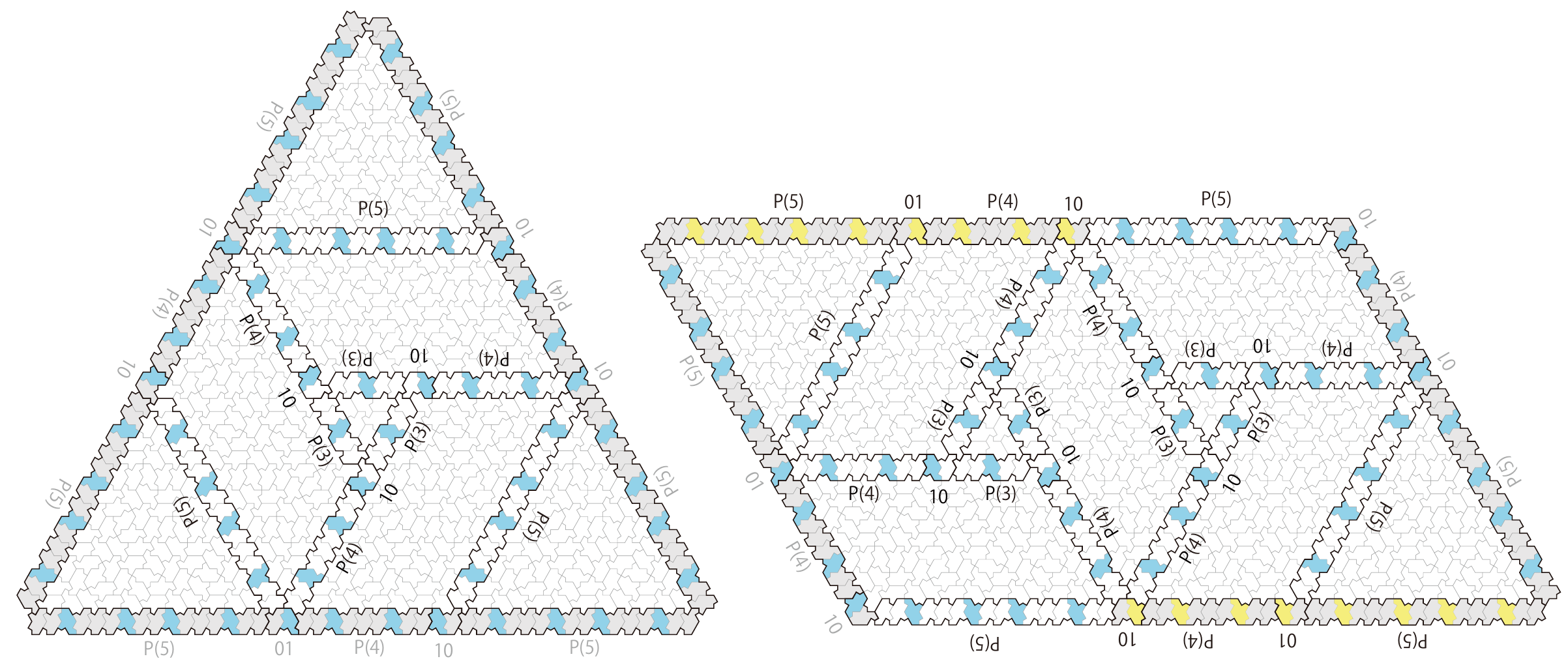}
\caption{Figure \ref{T} and \ref{P} for $n=3$\label{TPlevel3}}
\end{figure}

We show that Figure \ref{T} and \ref{P} complete our induction\footnote{One can start from $n=4$ to avoid using a single point tile $T_0$.}
for $n\ge 3$.
%more or less, a proof without words. 
The formulas (\ref{Cent1}) and (\ref{Cent2}) show that the Figure \ref{T} and \ref{P} combinatorially sound.
In the broken lined parts of $\P_n$, Lemma \ref{Pal} allows us to rotate $P(2n-1)$ and $P(2n-2)$
by the angle $\pi$.
Therefore it remains to show the consistency 
at places where several patch-tiles 
are meeting, i.e., we have to show that the geometric realization of $01$ and $10$ geometrically fit at the indicated places
in Figure \ref{T} and \ref{P}. We directly see that it is valid for $n=2,3$ in 
Figure~\ref{GH2} and \ref{TPlevel3}. 
%Further $\T_n,\P_n$ with $n=2,3$ satisfies the condition of 
%Figure~\ref{T} and \ref{P}. 
%See Figure \ref{TPlevel3} for the case $n=3$.
For $n\ge 4$, we claim 
that the local situation around $01$ or $10$ is exactly the same as in level $n-1$. 
Let us look at $\rotatebox[origin=c]{180}{10}$ in the bottom line of
$\P_n$. This position is surrounded by 
$$
\rotatebox[origin=c]{180}{$\P_{n-1}$},
\rotatebox[origin=c]{120}{$\P_{n-1}$},
\rotatebox[origin=c]{60}{$\P_{n-1}$} \quad \text{and} \quad    
\rotatebox[origin=c]{180}{$P(2n)$}
$$
and the shape is decided by
$$
\rotatebox[origin=c]{180}{$P(2n-2)$},
\rotatebox[origin=c]{120}{$P(2n-2)$},
\rotatebox[origin=c]{60}{$P(2n-2)$} \quad \text{and} \quad    
\rotatebox[origin=c]{180}{$P(2n-2)$}.
$$
By (\ref{Cent1}) and (\ref{Cent2}), it is surrounded by
$$
\rotatebox[origin=c]{180}{$P(2n-4)$},
\rotatebox[origin=c]{120}{$P(2n-4)$},
\rotatebox[origin=c]{60}{$P(2n-4)$} \quad \text{and} \quad
\rotatebox[origin=c]{180}{$P(2n-4)$}
$$
and $\rotatebox[origin=c]{180}{10}$ fits in this place in $\P_{n-1}$. Therefore
the consistency at this place is seen in the
induction assumption.
For the GSP 
$\rotatebox[origin=c]{180}{10}$ in the bottom of $\P_{n-1}$ which is 
located in the top right of $\P_{n}$, it is surrounded by 
$$
\P_{n-1}, \rotatebox[origin=c]{180}{$\T_{n-2}$},
\rotatebox[origin=c]{60}{$\P_{n-1}$}.
$$
Thus its shape is surrounded by
$$
\rotatebox[origin=c]{180}{$\P_{n-2}$},
\rotatebox[origin=c]{120}{$\P_{n-2}$},
\rotatebox[origin=c]{60}{$\P_{n-2}$},
\rotatebox[origin=c]{180}{$P(2n-4)$},
\rotatebox[origin=c]{180}{$\T_{n-2}$}
\quad \text{and} \quad
\rotatebox[origin=c]{60}{$\P_{n-1}$}
$$
and decided by
$$
\rotatebox[origin=c]{180}{$P(2n-4)$},
\rotatebox[origin=c]{120}{$P(2n-4)$},
\rotatebox[origin=c]{60}{$P(2n-4)$},
\rotatebox[origin=c]{180}{$P(2n-4)$},
\rotatebox[origin=c]{180}{$\T_{n-2}$}
\quad \text{and} \quad
\rotatebox[origin=c]{60}{$\P_{n-1}$}.
$$
By (\ref{Cent1}), (\ref{Cent2}), and Figures \ref{T} and \ref{P}, 
it is surrounded by
$$
\rotatebox[origin=c]{180}{$P(2n-6)$},
\rotatebox[origin=c]{120}{$P(2n-6)$},
\rotatebox[origin=c]{60}{$P(2n-6)$},
\rotatebox[origin=c]{180}{$P(2n-6)$},
\rotatebox[origin=c]{180}{$\T_{n-3}$}
\quad \text{and} \quad
\rotatebox[origin=c]{60}{$\P_{n-2}$}.
$$
We see $\rotatebox[origin=c]{180}{10}$ fits at the corresponding 
place in $\P_{n-1}$ by the
induction assumption.\footnote{Inspecting in detail, $\rotatebox[origin=c]{180}{$\T_{n-3}$}$ is irrelevant in the last statement.}
All occurrences of 
$01$ and $10$ in Figure \ref{T} and \ref{P} are checked in the same manner.
The claim is proved and our induction is completed.
\end{proof}

\section{Golden Ammann bar and Aperiodicity}
\label{GoldenAmmann}

If a tiling is invariant by a translation of a vector $v\in \R^2$, then $v$ is a period of the tiling. If any period of the tiling must be zero, 
we say that the tiling is non-periodic\footnote{See \cite[Section 1.3]{SMKGS:23_1} for different definitions of non-periodicity.}.
%In \cite{Gruenbaum-Shephard:87}, it is shown that 
%if there exists a non-zero period of a tiling by a given set of tiles
%in $\R^2$, then we can find a
%possibly different tiling by the same set of tiles 
%having two linearly independent periods, i.e., being lattice periodic. Therefore in our setting, we do not have to distinguish two concepts: the existence of a period and lattice periodicity.

\begin{figure}[h]
\begin{center}
\includegraphics[width=5cm]{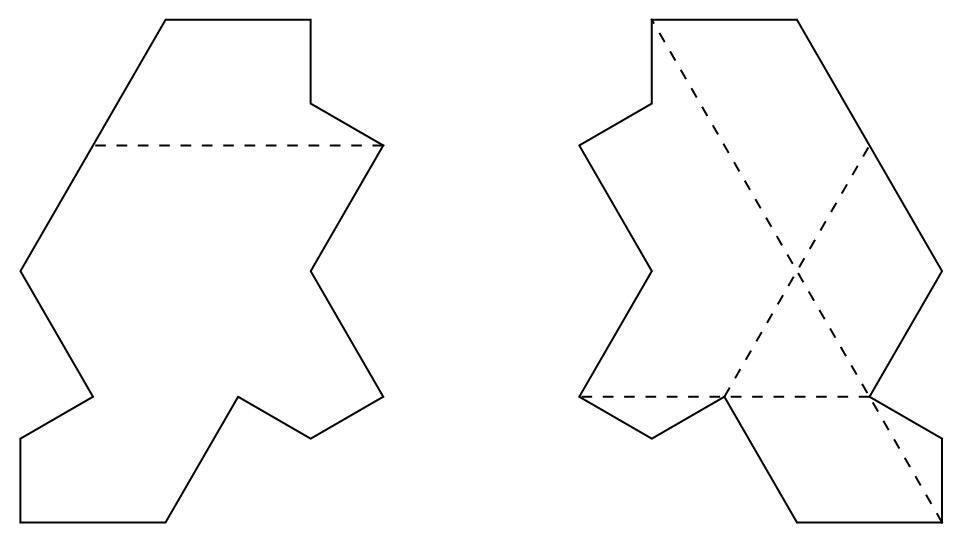}
\end{center}
\caption{Golden Ammann Barred tile\label{Amm0}}
\end{figure}

In this section, we will prove that
any tiling generated by Turtle is non-periodic. 
We introduce special markings in Figure \ref{Amm0}, which we call {\bf Golden Ammann Bar}s, in short {\bf GAB}s. 
We draw one dashed segment in 
the 
fore side, and three on the rear side.

Our strategy of proof is to show that there is 
an approximate hexagonal lattice structure 
of GAB within the Kagome Lattice
(Lemma \ref{AB},\ \ref{KagomeLattice} and \ref{Hexagon}).
Fixing two directions of GABs, 
the flipped tiles are exactly
located at the 
crossings of these GABs (Lemma \ref{Bij}). 
Since the ratio of lengths of 
GAB on the fore side and the rear is $1:4$, 
if GABs have natural frequency among Kagome tiling, 
then it must be irrational (Theorem \ref{AP}) by
statistical consideration of GABs measured by their
lengths. This implies that a
periodic tiling by Turtle
is impossible. 

Given a polygon that forms the boundary of a patch, 
the inner angle $\theta$ of a vertex is defined as usual. 
The ``complementary angle" of a vertex is defined to be $2\pi-\theta$, 
which is not the external angle $\pi-\theta$.

\begin{figure}[h]
\begin{center}
\includegraphics[width=8cm]{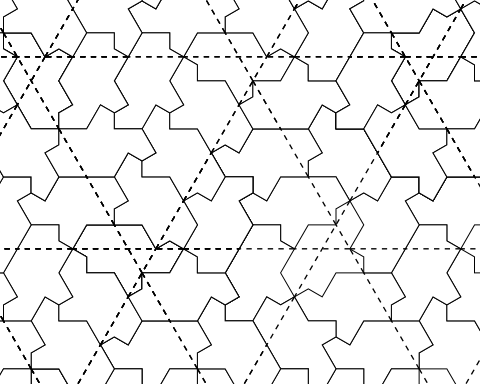}
\end{center}
\caption{Configuration by Golden Ammann Bars\label{AmmPatch0}}
\end{figure}

\begin{lem}
\label{AB}
For any tiling by Turtle, the
Ammann Bar in Figure \ref{Amm0} must continue
across the boundary to form a line
as in Figure \ref{AmmPatch0}.
\end{lem}

\begin{proof}
The set of inner angles of the
Turtle is $\{\pi/2,\ 2\pi/3,\ 4\pi/3,\ 3\pi/2,\ \pi\}$. Here the last
$\pi$ means the angle of a vertex on the edge.
Clearly, if a patch contains an acute complementary angle 
then it is impossible to 
extend to a tiling.

We observe that all endpoints of GAB in a Turtle
are located in the middle of an edge or a vertex whose angle or complementary angle is the right angle. Thus the outward
extension of this GAB must be covered by an edge or a right-angle vertex, i.e.,
the vertex of angle $2\pi/3$ and $4\pi/3$ does not help this covering.
%These simple observations are enough to see that any tiling by
%Turtle, each GAB must continue across the boundary as in Figure \ref{AmmPatch0}Therefore, 
At the endpoint of GAB, the possible angle configurations are 
$$\pi/2 +\pi/2+\pi/2+ \pi/2,\  
\pi/2+3\pi/2,\  \pi/2+\pi/2+\pi,\  \pi+\pi.$$ 

{\it Case 1.} For $\pi/2+3\pi/2$, there is no configuration which 
interrupt GAB at this outward extension, because of the edge length restriction, or the angles fit but the remaining configuration is impossible to continue, see Figure \ref{2Rights}.

\begin{figure}[h]
\begin{center}
\includegraphics[width=10cm]{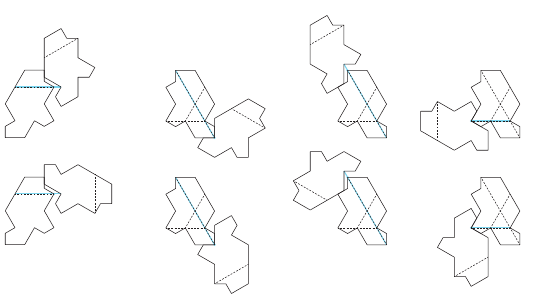}
\end{center}
\caption{Impossible to interrupt GAB at $\pi/2+3\pi/2$\label{2Rights}}
\end{figure}

{\it Case 2.}
For $\pi/2+\pi/2+\pi$ or $\pi+\pi$, the outward extension must be covered by a
vertex on the edge. 
To avoid an impossible configuration, 
this is possible only when the edge is the longest one. 
If it is not the beginning of another GAB, 
then the remaining configuration has an acute complementary angle.
If it is the beginning of another GAB, then the GAB continues straight there.

{\it Case 3.} For $\pi/2 +\pi/2+\pi/2+ \pi/2$, 
if GAB continues but does not go straight, then the tile having an outward 
GAB and the original tile form a configuration that has an acute complementary angle and this is out of this case. 
Therefore we only have to study the case that GAB has no outward connection. 
\begin{figure}[h]
\begin{center}
\includegraphics[width=10cm]{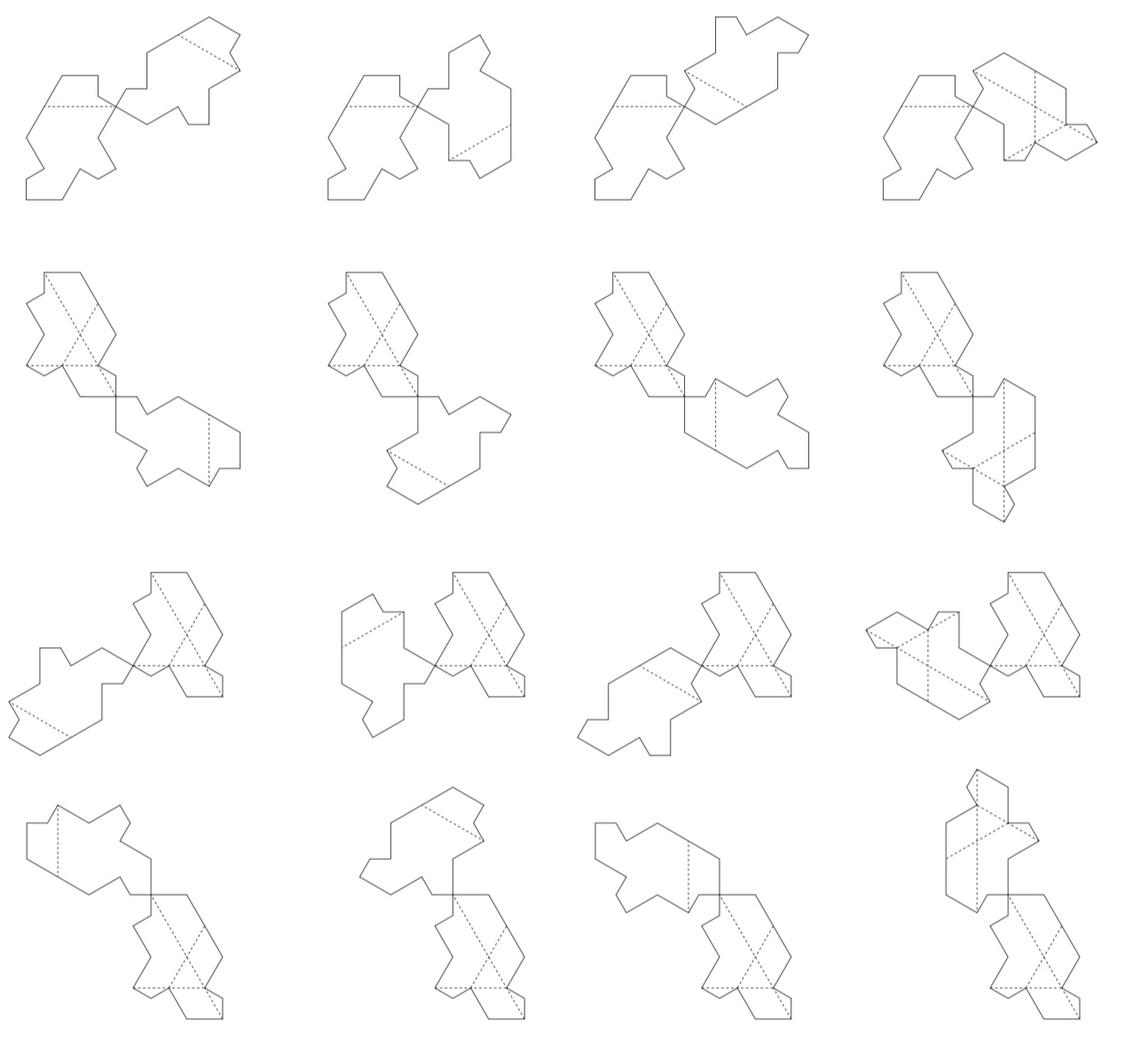}
\end{center}
\caption{Possible GAB disconnection\label{Disconnected}}
\end{figure}
There are $16$ such configurations as in Figure \ref{Disconnected}
that $\pi/2$ angle vertices meet 
but GAB does not continue. One can immediately confirm 
that none of the configurations extends to a tiling.
%For safety, we also 
%checked this property of GAB
%by listing all patches of minimum cardinality
%which contains an endpoint of GAB as an inner point.

%Such markings are originally introduced by R.~Ammann to construct a set of aperiodic tiles. Note that the situation around aperiodicity is very different. 
%The tiles by Ammann admits a periodic tiling without the marking bar, but this GAB is expected to be not necessary in the tiling of Tile $(\sqrt{3}:1)$. 
\end{proof}

\begin{rem}
\label{Ammann}
This type of markings was
originally introduced by R.~Ammann to construct 
a finite set of tiles to enforce hierarchical substitutive structure
 by additionally 
assuming that they must continue across the edges to a line, see
 \cite{AGS:92,Gruenbaum-Shephard:87,Akiyama:12}.
We should note that the role of Ammann bars in a Turtle is 
different because our GAB are ``dispensable" by Lemma \ref{AB}. 
We draw them only to show that the resulting tiling is non-periodic without knowing its substitutive structure.
\end{rem}

Our GAB serves supplementary information 
for the proof of the aperiodicity of Turtle. Hereafter we assume 
the edge length of the
small regular triangle formed by GABs of 
the flipped tile is $1/2$.

Let us fix two directions of Ammann bars and consider the set of crossing 
points. Here is an important observation:

\begin{lem}
\label{Bij}
There is a bijection between flipped tiles (the right of 
Figure \ref{Amm0}) and the set of intersection
of Ammann bars in the fixed two directions.
\end{lem} 

Figure \ref{AmmPatch0} may help our understanding of this statement.

\begin{proof}
Four segments
of length $1/2$ are emanating from
the crossing of two GABs.
It is impossible to cover this local 
configuration with the GABs of four non-flipped tiles.
There exists at least one flipped tile. 
However, once we use a
flipped tile to cover this crossing, three or four segments are covered
among the four segments.
Therefore this crossing point must be covered in one of two ways.
It is either 
covered by a single flipped tile, or by exactly one flipped tile and 
one non-flipped tile. As a result, 
each crossing is contained in exactly one flipped tile.
This gives a map from a crossing to a flipped tile. Since two directions
are fixed, we can recover uniquely the crossing point from the flipped tile.
\end{proof}

% when all crossings are found 
%in small regular triangles consisting of GABs.
%\footnote{
%Some crossings may lie outside this small GAB triangle when the GAB appears more often. See the complementary GABs in Figure \ref{AmmPatch}.}. 

\begin{figure}[h]
\begin{center}
\includegraphics[width=5cm]{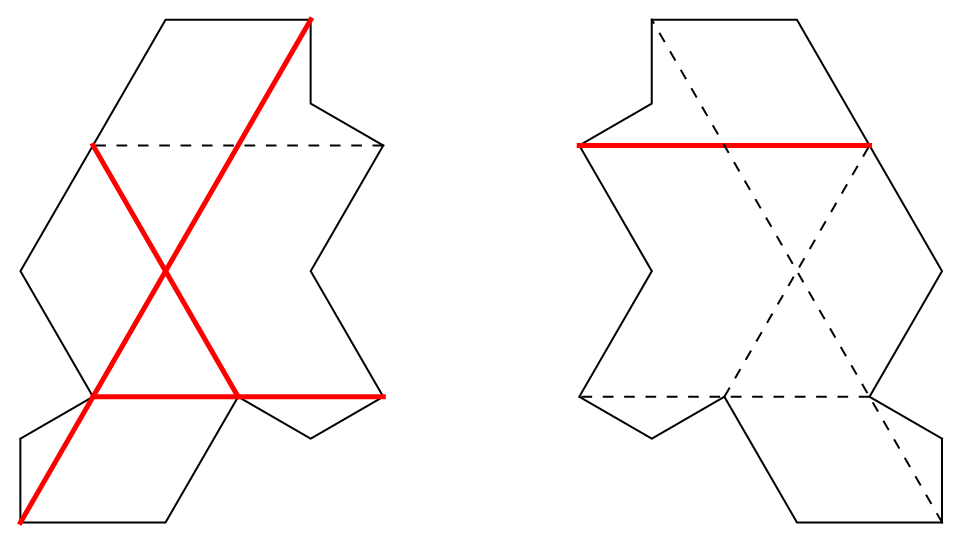}
\end{center}
\caption{Complementary Golden Ammann Bars\label{Amm}}
\end{figure}

One can also define complementary markings as in Figure \ref{Amm}; three red segments on the fore-side and one segment in the rear.

A generalized GAB is either a GAB or a complementary GAB. 

\begin{figure}[h]
\begin{center}
\includegraphics[width=8cm]{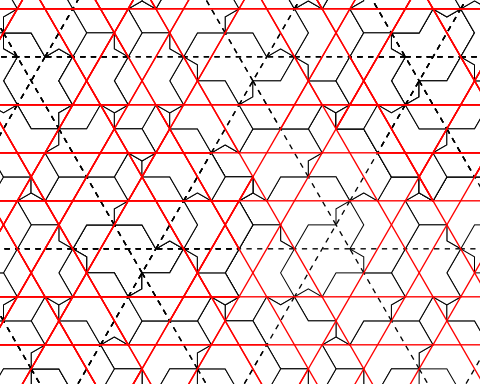}
\end{center}
\caption{All Golden Ammann Bars\label{AmmPatch}}
\end{figure}

\begin{figure}[h]
\begin{center}
\includegraphics[width=8cm]{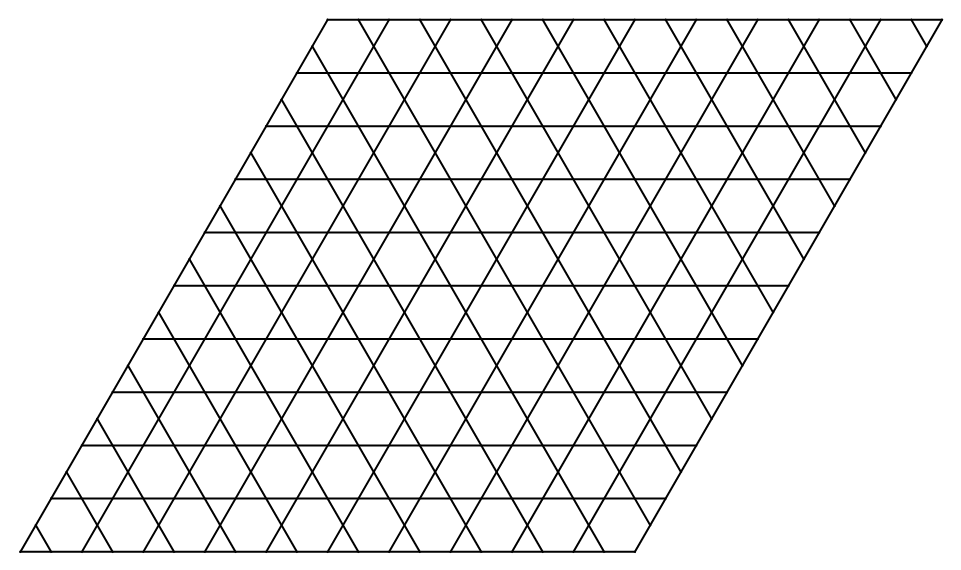}
\end{center}
\caption{Kagome parallelogram $K(10)$\label{Kagome}}
\end{figure}

\begin{lem}
\label{KagomeLattice}
The set of generalized GABs  
must form a `Kagome' tiling (also called trihexagonal tiling) 
as in Figure \ref{AmmPatch} or \ref{Kagome}, 
one of the 2-uniform tilings whose signature is
$3636$\footnote{This generalized GAB and Figure \ref{AmmPatch}
were
observed in \cite{Reitebuch:23}.}.
\end{lem}

\begin{proof}
By Lemma \ref{AB}, the generalized Ammann 
Bars must continue to straight lines in a tiling.
For every tiling by Turtle, one can show  that all
adjacent parallel generalized GABs are separated by distance 
$\sqrt{3}/2$. For example, Figure \ref{KAB} shows
all possible ways to fill the black spot, which forces three parallel
generalized GABs. 
One can easily do this case analysis for all 6 directions.
\begin{figure}[h]
\begin{center}
\includegraphics[width=10cm]{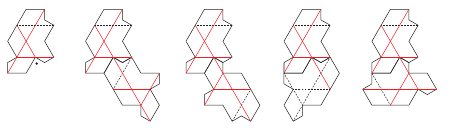}
\end{center}
\caption{Kagome Structure\label{KAB}}
\end{figure}
Having this in mind, to form a small regular triangle of edge length $1/2$
on a tile, the generalized GAB must form the Kagome tiling.
\end{proof}

We consider the Kagome GAB parallelogram $K(n)$ in Figure \ref{Kagome}
consisting of $(n+1)$ segments in the horizontal direction, $(n+1)$ segments of
slope $\sqrt{3}$ and $2n$ segments of slope $-\sqrt{3}$.
By assumption, 
the small regular triangles in the Kagome pattern $K(n)$ has edge length $1/2$.

We associate a
hexagonal coordinate $\langle x,y \rangle\in \R^2$ to the 
Kagome parallelogram by 
$\langle x,y\rangle:=x(1,0)+y(1/2,\sqrt{3}/2)$. 
Segments
$H_i\ (i=0,1,\dots, n)$ connect 
$\langle 0,i\rangle$ to $\langle n,i\rangle$ and 
segments $L_j\ (j=0,1,\dots, n)$ of $\sqrt{3}$ slope
connect $\langle j,0\rangle$ to $\langle
j,n\rangle$. Finally
$M_k\ (k=0,\dots,2n-1)$ are segments of $-\sqrt{3}$ slope 
 which connect
$\langle 1/2+k,0\rangle$ and $\langle 0,k+1/2\rangle$ for $k\in \{0,1,\dots,n-1\}$, 
and $\langle n,k-n+1/2\rangle$ and $\langle k-n+1/2,n\rangle$ for $k\in \{n,n+1,\dots, 2n-1\}$.

Let us fix a tiling by Turtle such that a vertex of its Kagome tiling is 
at the origin.
Assume that $(a_i)_{i=0}^{\infty}$, $(b_j)_{j=0}^{\infty}$, $(c_k)_{k=0}^{\infty}$
are three increasing sequences, $a_i,b_j$ are non negative integers and $c_k$
is a half-integer so that horizontal
GABs go through 
$\langle 0,a_i\rangle$, GABs of slope $\sqrt{3}$
pass $\langle b_j, 0\rangle$ and the ones of slope $-\sqrt{3}$ 
pass $\langle c_k,0\rangle$.
Choose $h(n),\ell(n),m(n)$ so
that
$H_{a_i}\ (i=0,\dots, h(n))$, 
 $L_{b_j}\ (j=0,\dots, \ell(n))$ and $M_{c_k}\ (k=0,\dots, m(n))$ form 
the set of GABs in $K(n)$.
We may assume that $h(n)$ and $\ell(n)$ are 
positive. Indeed since one may take the mirror image
of the tiling, we may assume that there are infinitely many GABs in at least two directions, say in $H_i$ and $L_j$ directions.
%Since $H_{a_i}$ and $L_{b_j}$ intersect at $\langle a_i,b_j
%\rangle$, there exists a unique 
%flipped tile whose GABs cover three or four segments of length $1/2$
%emanating from 
%$\langle a_j,b_j\rangle$. 
%Therefore by the definition of GAB, 
By Lemma \ref{Bij}, there exists a unique $k$ such that 
$c_k-a_i-b_j$ is equal to $\pm 1/2$ 
and $M_{c_k}$ is passing the flipped tile 
where $H_{a_i}$ and $L_{b_j}$ have 
the crossing.
There is a natural ordering that
if  $i\le i'$, $j\le j'$ and $c_{k'}-a_{i'}-b_{j'}=\pm 1/2$ then $k\le k'$.
Symmetric discussion holds 
for the intersection of $H_{a_i}$ and $M_{c_k}$, as well as for 
the intersection of $L_{b_j}$ and $M_{c_k}$. 
Next lemma shows an approximate 
hexagonal lattice structure of the flipped tiles.

\begin{lem}
\label{Hexagon}
We obtain a relation $k=i+j$, i.e., a formula 
\begin{equation}
\label{Hex}
c_{i+j}-a_i-b_j = \pm 1/2.
\end{equation}
Consequently we have $m(n)=h(n)+\ell(n)$. 
Note that the sign $\pm$ in (\ref{Hex}) depends on $i$ and $j$.
\end{lem}

\begin{proof}
Let $c_k-a_0-b_0=\pm 1/2$. The above
symmetric logic gives $c_0-a_{\ell}-b_0=\pm 1/2$ with a unique $\ell$. By the natural ordering, we have $k=\ell=0$. 
If $c_k-a_i-b_j=\pm 1/2$ and $k\neq i+j>0$, then 
take the minimum $i+j$ with this property. 
Since $k<i+j$ gives a contradiction to the uniqueness of $k$, we have $k>i+j$. Using the above symmetric logic, there exists a unique $\ell$ with 
$c_{i+j}-a_{\ell}-b_j=\pm 1/2$. By the induction hypothesis, we see $\ell\ge i$.
By the natural ordering, we infer $\ell\le i$ and thus $\ell=i$. 
This gives a contradiction to the uniqueness of $k$. Lemma \ref{Hexagon} is proved.\end{proof}

%Therefore
%\begin{align*}
%c_{i+j}-a_i-b_j&=\pm 1\\
%c_{i+j}-a_{i+1}-b_{j-1}&=\pm 1
%\end{align*}
%gives
%$$
%a_{i+1}-a_i -(b_j-b_{j-1})\in \{-2,0,2\}.
%$$
%Similar discussion yields
%\begin{align*}
%&c_{i+1}-c_{i}-(a_{j+1}-a_j) \in \{-2,0,2\}\\
%&c_{i+1}-c_i-(b_{j+1}-b_j) \in \{-2,0,2\}.
%\end{align*}
%This shows that dashed GABs must have a structure close to a hexagonal
%lattice.
%The gap between adjacent GABs of the same orientation is close to a constant 
%function, see Figure \ref{AmmPatch0}.\footnote{These results suggest that the set of all gaps $a_{i+1}-a_i$,$b_{i+1}-b_i$ and
%$c_{i+1}-c_i$ may be small in cardinality, 
%but this kind of fact is not necessary to show non-periodicity of GAB.}
%The situation is the same for the complemetary GAB as in Figure \ref{AmmPatch}.
Let $\N$ be the set of non-negative integers and 
$U$ be a subset of $\N$. If
$$
\delta(U):=
\lim_{n\to \infty}
\frac 1{n+1} \mathrm{Card}(U \cap [0,n])
$$
exists then $\delta(U)$ is called the natural density of $U$. We prove the following lemma
\begin{lem}
\label{Abel}
The natural density $\delta(U)$ exists if and only if 
$$
\lim_{n\to \infty} \frac 1{n^2} \sum_{j\in U\cap [0,n]} j
$$
exists. In this case, the last limit is equal to $\delta(U)/2$.
\end{lem}

\begin{proof}
%We show it for the convenience of readers.
Let $\chi_U$ be the indicator function of $U$. 
Then we see
$$
\delta(U)=\lim_{n\to \infty} \frac 1{n+1}
\sum_{j=0}^n \chi_U(j).
$$
Assume that $S_n=\sum_{j=0}^n \chi_U(j)j=n^2 q/2 +o(n^2)$. Then
\begin{align*}
\sum_{j=0}^n \chi_U(j)
%&=\chi_U(0)+ \sum_{j=1}^n \chi_U(j)\\
&=\chi_U(0)+ \sum_{j=1}^n \frac{S_j-S_{j-1}}j\\
&=\frac {S_n}n+ \chi_U(0)+\sum_{j=1}^{n-1} \frac {S_j}{j(j+1)}\\
&=\frac {nq}2 + o(n)
+\chi_U(0)+ \sum_{j=1}^{n-1} \left(\frac q2\left(1-\frac 1{j+1}\right)+o(1)\right)\\
&=\frac {nq}2 + \frac {(n-1)q}2 - \frac {q\log n}2 + o(n)\\
&=nq + o(n).
\end{align*}
For the converse, assume that $T_n=\sum_{j=0}^n \chi_U(j)=n q +o(n)$.
Then we have
\begin{align*}
\sum_{j=0}^n \chi_U(j)j &= \sum_{j=1}^n (T_j-T_{j-1})j\\
&=n T_n -\sum_{j=1}^{n-1} T_j\\
&=n^2 q + o(n^2) -\sum_{j=1}^{n-1} (jq + o(j))\\
&=\frac {n^2q}2+o(n^2).
\end{align*}
\end{proof}

We are ready to state the second main result of this paper.

\begin{thm}
\label{AP}
For a tiling by Turtle, Golden Ammann Bars give a sub-configuration as in 
Figure \ref{AmmPatch0} of Kagome tiling (Figure \ref{Kagome}). 
If Kagome lines become GAB with frequency $q$ in one direction, then it has 
the same frequency $q$ in the
other two directions. In this case, the value $q$
must be equal to one of $\frac {5\pm \sqrt{5}}{10}$. 
Consequently, each tiling by Turtle 
is non-periodic.
\end{thm}

\begin{proof}
The first statement summarizes Lemma \ref{AB} and \ref{KagomeLattice}.
Assume that if GABs in the horizontal direction has the frequency $q$
among the horizontal lines of Kagome tiling, 
i.e., there exists $q\in [0,1]$ such that
\begin{equation}
\label{Freq0}
\lim_{n\to \infty} \frac {h(n)}{n} = q.
\end{equation}
Shifting the tiling, we may assume that $a_0=b_0=0$.
By Lemma \ref{Hexagon}, we have
\begin{align*}
&c_{j}-a_j\in [-1/2,1/2],\\
&c_{j}-b_j\in [-1/2,1/2],\\
&a_j-b_j\in [-1,1].
\end{align*}
%Taking $n\to \infty$, we may extend 
%$(a_i),(b_j),(c_k)$ to 
%infinite sequences having non-negative integers
%indexes $i,j,k$.
Setting $h=h(n)$, from $a_h\le n < a_{h+1}$, we have 
$$
b_h\le n+1,\quad b_{h+1}\ge n.
$$
This implies
%$$
%b_{h-t+1}\le n-t+2,\quad b_{h+t}\ge n+t-1
%$$
%holds for $t=1,2,\dots$. Thus we have
$
h(n)-1 \le \ell(n)\le h(n)+1
$
and
\begin{equation}
\label{Freq1}
\lim_{n\to \infty} \frac {\ell(n)}{n} = q.
\end{equation}
Since
\begin{align*} 
|c_k-a_k| &\le 1/2 & k\le h(n)\\
|c_k-(a_{h(n)}+b_{k-h(n)})|&\le 1/2 & h(n)<k\le m(n),
\end{align*}
$c_k$ and $a_k$ are in one to one correspondence in $k\le h(n)$ 
and $c_k$ and $b_{k-h(n)}$ are one to one in $k>h(n)$, because of
Lemma \ref{Hexagon} and the natural ordering. 
Thus the integer sequence 
$(c_k+1/2)$ inherits the frequency of $(a_i)$ and $(b_j)$ and we see
\begin{equation}
\label{Freq2}
\lim_{n\to \infty} \frac {m(n)}{n} = q.
\end{equation}
Thus the second assertion is proved.

The area of a
Turtle is $13+1/3=40/3$ times the area of the small regular triangle in Kagome tiling and $K(n)$ consists of $8n^2$ small regular triangles. 
This shows that
the minimum number of Turtles which covers $K(n)$ is
\begin{equation}
\label{Count0}
\frac{3}{5} n^2 + O(n).
\end{equation}
Since there are $O(n)$ tiles which intersect the outermost parallelogram of $K(n)$, the number of Turtles lie strictly within $K(n)$ is also 
(\ref{Count0}).
Later we shall implicitly use this fact that the number of Turtles is insensitive to the way we count them, up to this error term.

We will
compute the sum of all lengths of GABs in two ways.
By (\ref{Freq0}), the sum of length of $H_{a_i}$ are
$$
n((n+1)q +o(n))= n^2 q +o(n^2).
$$
%Note that the distance between two parallel lines in $K(n)$ is $2$ by definition. 
The same is valid for $L_{b_j}$ by (\ref{Freq1}).
The sum of length of $M_{c_k}$ is divided into two parts:
$$
\sum_{c_k<n} c_k + \sum_{c_k\ge n} 
\left(2n-c_k\right).
$$
By using (\ref{Freq2}) 
and Lemma \ref{Abel} with $U=\{c_k+1/2\ |\ c_k<n\}$, 
we have
$$
\sum_{c_k<n} \left(c_k+\frac 12\right) = \frac{n^2 q}2 + o(n^2).
$$
Similar computation shows
$$
\sum_{c_k\ge n} \left(2n-c_k-\frac 12\right)= 2n^2q- \left(\frac{(2n)^2 q}2
-\frac{n^2 q}2\right) + o(n^2)=\frac{n^2 q}2 + o(n^2).
$$
The contributions of $\pm 1/2$ in these two formulas are $O(n)$. 
Therefore we see that the total length of GABs is
\begin{equation}
\label{Count1}
3 n^2 q + o(n^2).
\end{equation}
By Lemma \ref{Bij}, 
the crossing of $H_{a_i}$ and $L_{b_j}$ uniquely corresponds to 
a flipped tile. Thus we find $n^2 q^2+O(n)$ flipped tiles in $K(n)$. 
%Consider the parallelogram tiling $P(n)$
%generated by $H_{a_i}$ and $L_{b_j}$, a sub-configuration of $K(n)$. 
%Every crossing of $P(n)$ is already covered 
%by the flipped tiles and the remainder has no crossing of GABs in $P(n)$. 
%This means in $K(n)$ we can not find a small triangle after the
%removal of
%the above flipped tiles. Therefore the above 
%$n^2 q^2+O(n)$ flipped tiles are the totality of flipped tiles in $K(n)$.
The length of GAB on the fore-sided
tile is $1$ while the length is $4$ in the rear side. Using (\ref{Count0}),
the total length of GAB in $K(n)$ is computed in a different way:
\begin{equation}
\label{Count2}
\left(\frac{3}{5} n^2 -n^2 q^2 +o(n^2)\right) + 4 \left(n^2 q^2 +o(n^2)\right)
= \left(\frac {3}5+3 q^2 \right)n^2 + o(n^2).
\end{equation}
Comparing (\ref{Count1}) and (\ref{Count2}) as $n\to \infty$, 
we obtain
$$
q^2-q+\frac 15=0.
$$
Thus
$$
q_1=\frac{5-\sqrt{5}}{10}=\frac 1{1+\tau^2}\approx 0.276393, \quad 
q_2=\frac{5+\sqrt{5}}{10}=\frac {\tau^2}{1+\tau^2}\approx 0.723607
$$
are the possible two values of frequency $q$. Clearly 
both of them are irrational. 

If there exists a period $v\neq 0$ of a tiling by Turtle, then 
every tile is sent to the tile of the same orientation. 
Thus the set of GABs is invariant by this translation. Rotating the tiling 
if necessary, we may assume that the horizontal frequency $q$ 
must exist and it is rational, we obtain a contradiction. 
\end{proof}

Note that two values correspond to the frequency of GAB and that of 
complementary GAB, and $q_2/q_1=\tau^2$. The Golden Hex tiling in the previous section has the GAB frequency $(5- \sqrt{5})/10$, which is easily shown by the existence of arbitrary long $\GSE$'s.

\begin{rem} 
After the submission to ArXiv, 
we are informed that a different 
proof of aperiodicity using GAB was released 
several days prior to our post, see \cite{MathBlock}. Here is the main difference.
The proof in \cite{MathBlock} relies on an assumption that the tiling by Turtle must have an underlying $[3.4.6.4]$ Laves tiling.
A proof of this assumption is found in \cite{SMKGS:23_1} 
for Smith Hat ($b=1/\sqrt{3}$), but 
the one for Turtle is postponed in \cite{MathBlock}. In contrast, 
our proof is self-contained.
It seems Lemma \ref{Hexagon} plays the role of the assumption. 
Ideas of two proofs are also different; we computed the total
length of GAB, while \cite{MathBlock} computed the number of essential 
crossings of GABs and that of complementary GABs
%using angle consideration 
in the Laves tiling.
%Another small difference exists in the definitions of non-periodicity: we showed that any period must be zero, while \cite{MathBlock} denied 
%the existence of two linearly independent periods.
\end{rem}

\section{Future works}
\label{Open}

Two substitutions are proposed in \cite{SMKGS:23_1} based on the extensive search of the corona shapes. 
One of them is based on a patch-tiles H7 and H8, 
and it is of interest to study the existence of the
related cut and project scheme of the limit H7H8 tiling.
We checked the pure discreteness of limit H7H8 tiling dynamics 
by the algorithm in \cite{Akiyama-Lee:10}. 
This shows that 
there exists a $2\times 2$ cut and project scheme (c.f.\ \cite{BaakeGaehlerSadun:23}).
We originally found the Golden Hex substitution structure
in the subdivision of the associated
cut and project window (c.f.\ \cite{Baake-Grimm:13,Akiyama-Barge-Berthe-Lee-Siegel}). We shall discuss this relationship elsewhere. 

From our Golden Hex substitution whose consistency is legitimated in the first section, we expect that the combinatorial substitution rule of H7 and H8 in \cite{SMKGS:23_1} is realized as a concrete geometric substitution. 
It is plausible that our Golden Hex tiling is MLD with H7H8 tiling
before taking their limits.
Our preliminary experiments show that it seems to be the case 
and therefore this Golden Hex substitutive structure may be 
automatically enforced for all tilings by Turtle.

Apart from Smith Hat in \cite{SMKGS:23_1}, 
several other aperiodic monotiles are studied
in \cite{SocolarTaylor:10, MampustiWhittaker:20, 
WaltonWhittaker:23,
SMKGS:23_2} under 
different conditions. 
The idea of the proofs of aperiodicity is 
to enforce some large structure in the resulting tilings.
Our proof using the
statistical property of GAB seems 
to be new and we expect some further applications.
\medskip

{\bf Acknowledgments.}
We would like to thank C.~Kaplan, F.~G\"ahler and M.~Baake
for helpful comments on the earlier version of the paper.
We are deeply indebted to the detailed comments from the referees. 
Due to this criticism, we could drastically improve the readability 
of the construction of Golden Hex substitution as well as 
the exposition in \S \ref{GoldenAmmann}. 
%The simpler proof of (\ref{Freq1}) 
%is also due to the referee. 
Thanks are also due to the study group for aperiodic tiles in Univ.~Tsukuba, 
in particular P.~Sonngam, P.~Wongpinta, Y-L.~Xu and K.~Ito who 
tirelessly checked the consistency of 27 vertex atlases of the
original 
substitution rule in \S \ref{Appe}, and gave us an idea for possible 
improvement.
\medskip

Data openly available in a public repository: https://arxiv.org/abs/2307.12322

\section{Appendix}
\label{Appe}

\begin{figure}[h]
\begin{center}
\includegraphics[width=10cm]{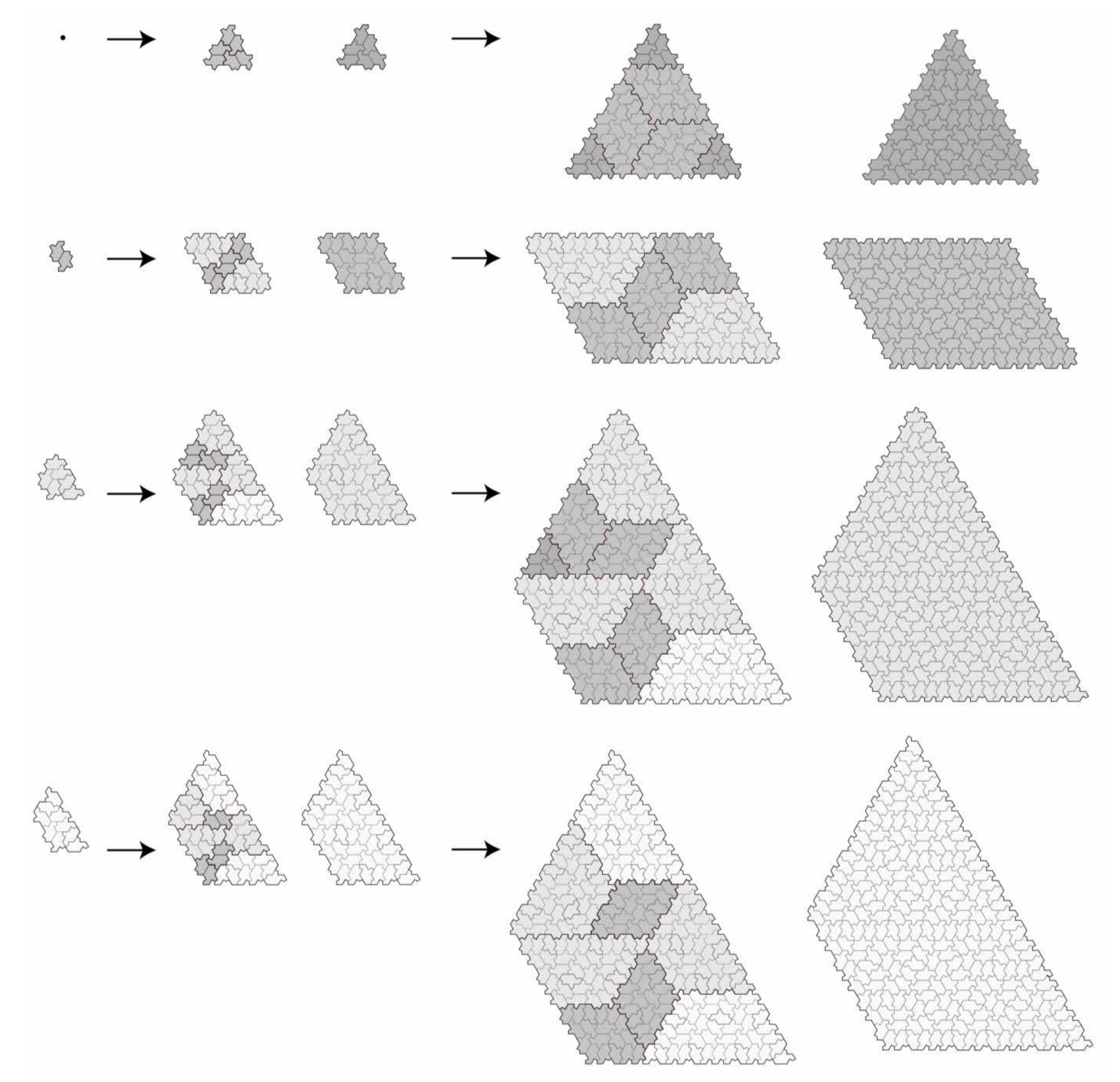}
\end{center}
\caption{$(\T_n,\PD_n,\ZD_n,\ZT_n)_{n=0,1,2,3}$\label{GH}}
\end{figure}

Our substitution in \S \ref{GoldenHex} became 
much simpler than the previous one (arXiv:2307.12322 ver 4).
We defined four growing sequences $(\T_n,\PD_n,\ZD_n,\ZT_n)_{n=0,1,\dots}$
of patch-tiles generated by
Turtle, whose first several terms are depicted in Figure~\ref{GH}.

\begin{figure}[h]
\begin{center}
\includegraphics[width=12cm]{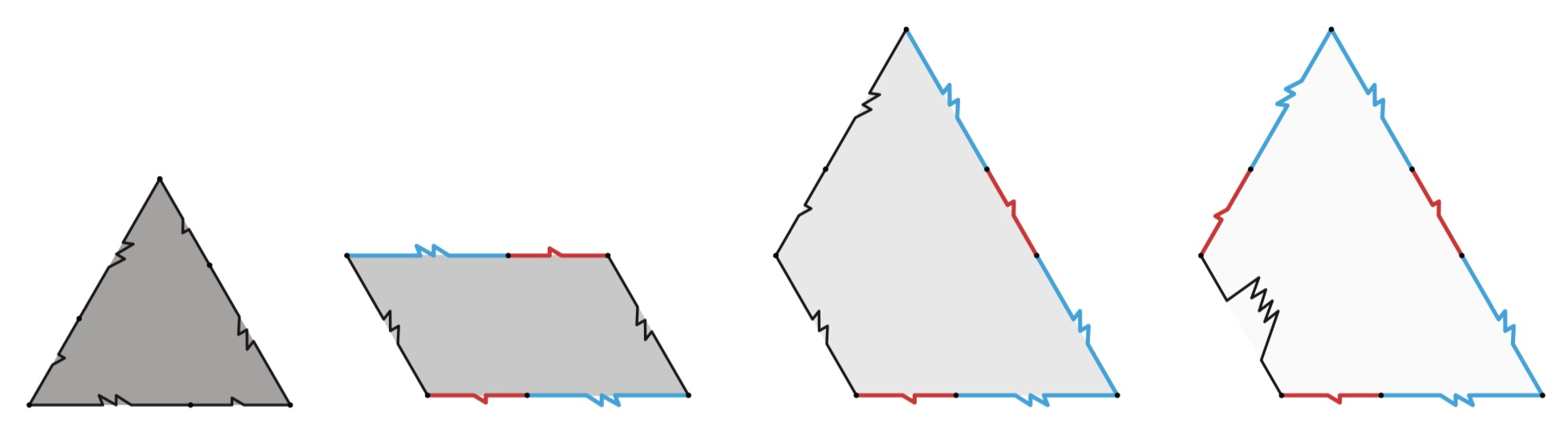}
\end{center}
\begin{center}
\includegraphics[width=12cm]{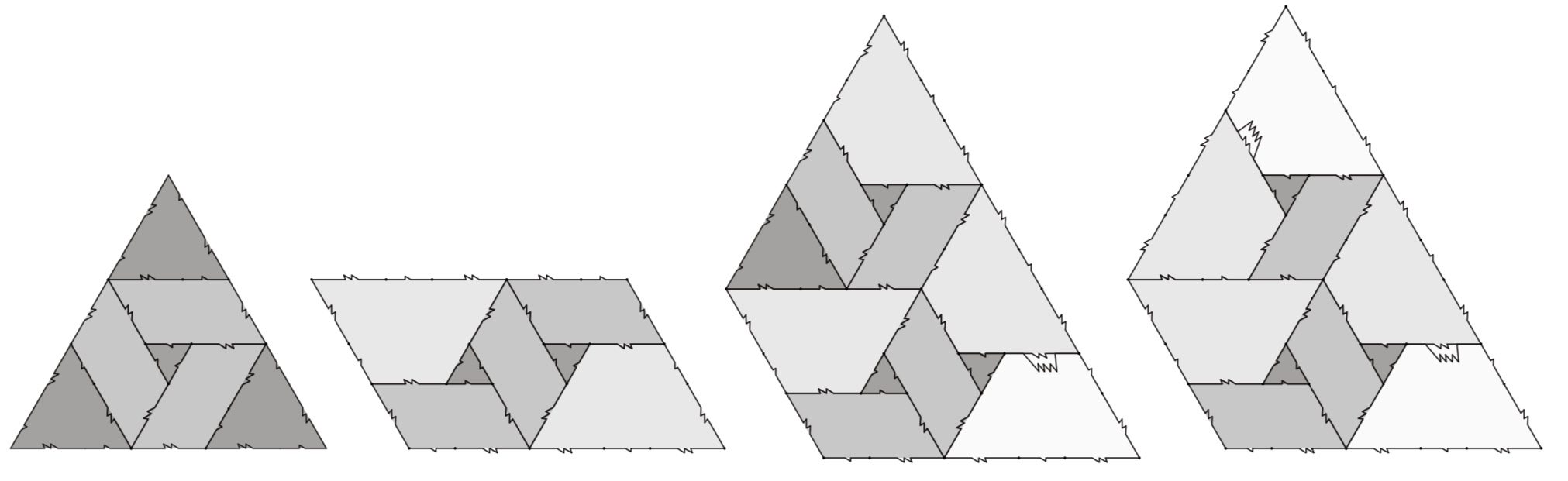}
\end{center}
\caption{Total Substitution \label{ES}}
\end{figure}

A vertex atlas is a patch in a given tiling 
that shares a common point (a vertex) 
on the boundary and the common point is an inner point 
of the union of the patch, having minimal cardinality with this property.
We collected 27 vertex atlases of patch-tiles of level 4
to {\bf define} the substitution and 
check that they give rise to the total substitution
in Figure \ref{ES} including their boundaries keeping its consistency.
Their boundaries are surrounded by GSPs, or receive GSPs, 
as indicated by notches and dents.
In the notation of \S \ref{GoldenHex}, one notched segment corresponds to $P(2n-1)$ and two notched one corresponds to $P(2n)10$. 
If the right end of $P(2n-1)$ finishes in the middle of the edge, then
we have to insert $01$ after $P(2n-1)$.
The four notched segment also 
corresponds to $P(2n)10$ but the left end is longer than the two notched one.
The shape looks ready to
receive $01P(2n)10$, but this does not happen in this substitution. The left
place for $01$ is occupied by a GSP of a different direction.
We see this fact by checking the 
27 vertex atlases\footnote{One can also check that
$1p_{2n}1$ is a forbidden word of the sturmian word of slope $(5-\sqrt{5})/10$.}.

The referee wrote us that this construction can not be called `simple'. 
This objection may be correct that it is pretty 
heavy to check and requires many precise drawings.

\begin{figure}[h]
\begin{center}
\includegraphics[width=12cm]{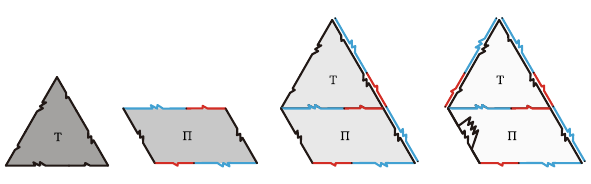}
\end{center}
\caption{$\P_n=PD_n$ \label{TS}}
\end{figure}

We claim that the set of tilings $\mathcal{C}$
generated by $(\T_n,\PD_n,\ZD_n,\ZT_n)_{n=0,1,\dots}$ is the same as the ones by $(\T_n,\P_n)_{n=0,1,\dots}$, which is denoted by $\mathcal{B}$.
By induction, we see that $\P_n=\PD_n$ by introducing subdivisions of $\ZD_n$ 
and $\ZT_n$ as in Figure~\ref{TS}. Here 
Lemma \ref{Pal} guarantees the 
$\pi$-rotation of GSPs on the boundary.
This justifies the abusive usage of 
the same symbol $\T_n$ in two substitutions.
Thus every patch of the tiling by $(\T_n,\P_n)_{n=0,1,\dots}$ must appear in
a tiling in $\mathcal{C}$. The total substitution in Figure \ref{ES} is
primitive, i.e., every patch appears as 
a subpatch of all $\T_n, \PD_n, \ZD_n, \ZT_n$ for some $n$. 
Since $\T_n$ is a patch of the tiling in 
$\mathcal{B}$ as well, 
every patch of the tiling by  
$(\T_n,\PD_n,\ZD_n,\ZT_n)_{n=0,1,\dots}$ appears in a
tiling of $\mathcal{B}$. The claim is proved.

By this discussion, the sequence 
$(\T_n,\PD_n,\ZD_n,\ZT_n)_{n=0,1,\dots}$ is reconfirmed to be well-defined through Theorem \ref{GHSTiling} for $(\T_n,\P_n)_{n=0,1,\dots}$.
If we include GSPs as patch-tiles, then
we lose primitivity of the substitution, because GSPs do not contain 
large balls.
It is noteworthy that the total substitution in Figure \ref{ES} is primitive and it did not treat GSPs as patch-tiles, while the  
rule in Figure \ref{T} and \ref{P} is not primitive. Since 
statistical and ergodic
properties are easier for primitive substitution 
than non-primitive one, 
the sequence
$(\T_n,\PD_n,\ZD_n,\ZT_n)_{n=0,1,\dots}$ is of independent interest.

%\bibliographystyle{amsplain}
%\bibliography{../reflist}

\begin{thebibliography}{10}

\bibitem{Akiyama:12}
S.~Akiyama, \emph{A note on aperiodic {A}mmann tiles}, Discrete Comput. Geom.
  \textbf{48} (2012), no.~3, 702--710.

\bibitem{Akiyama-Barge-Berthe-Lee-Siegel}
S.~Akiyama, M.~Barge, V.~Berth\'e, J.-Y. Lee, and A.~Siegel, \emph{On the
  {P}isot {S}ubstitution {C}onjecture}, Mathematics of Aperiodic Order,
  Progress in Mathematics, vol. 309, Birkh\"auser, Basel, 2015, pp.~33--72.

\bibitem{Akiyama-Lee:10}
S.~Akiyama and J.-Y. Lee, \emph{Algorithm for determining pure pointedness of
  self-affine tilings}, Adv. Math. \textbf{226} (2011), no.~4, 2855--2883.

\bibitem{AGS:92}
R.~Ammann, B.~Gr{\"u}nbaum, and G.~C. Shephard, \emph{Aperiodic tiles},
  Discrete Comput. Geom. \textbf{8} (1992), no.~1, 1--25.

\bibitem{BaakeGaehlerSadun:23}
M.~Baake, F.~G\"ahler, and L.~Sadun, \emph{Dynamics and topology of the hat
  family of tilings}, Ar{X}iv:2305.05639.

\bibitem{Baake-Grimm:13}
M.~Baake and U.~Grimm, \emph{Aperiodic {O}rder. {V}ol. 1}, Encyclopedia of
  Mathematics and its Applications, vol. 149, Cambridge University Press,
  Cambridge, 2013.

\bibitem{Fogg:02}
N.~{Pytheas Fogg}, {\em Substitutions in Dynamics, Arithmetics and
  Combinatorics}, vol.~1794 of {\em Lecture Notes in Mathematics}.
\newblock Springer-Verlag, 2002.

\bibitem{Frank-Sadun:14}
N. P. Frank and L. Sadun. \emph{Fusion: a general framework for hierarchical tilings of $\R^d$},
Geom. Dedicata \textbf{171} (2014), 149-186.

\bibitem{Gruenbaum-Shephard:87}
B.~Gr{\"u}nbaum and G.~C. Shephard, \emph{Tilings and patterns}, W. H. Freeman
  and Company, New York, 1987.

\bibitem{Lothaire:02}
M.~Lothaire, \emph{Algebraic combinatorics on words}, Encyclopedia of
  Mathematics and its Applications, vol.~90, Cambridge University Press,
  Cambridge, 2002.

\bibitem{MathBlock}
\emph{The Turtle prototile is not periodic, a simple proof},
\url{https://archive.li/rKD0U}

\bibitem{MampustiWhittaker:20}
M.~Mampusti and M.~F. Whittaker, \emph{An aperiodic monotile that forces
  nonperiodicity through dendrites}, Bull. Lond. Math. Soc. \textbf{52} (2020),
  no.~5, 942--959.

\bibitem{Reitebuch:23}
U.~Reitebuch, \emph{Direct construction of aperiodic tilings with the hat
  monotile}, Ar{X}iv:2306.06512.

\bibitem{SMKGS:23_1}
D.~Smith, J.~S. Myers, C.~S. Kaplan, and C.~Goodman-Strauss, \emph{An aperiodic
  monotile}, Combinatorial Theory \textbf{4.1} (July 2024).

\bibitem{SMKGS:23_2}
\bysame, \emph{A chiral aperiodic monotile}, arXiv:2305.17743.

\bibitem{SocolarTaylor:10}
J.~E.~S. Socolar and J.~M. Taylor, \emph{An aperiodic hexagonal tile}, Journal
  of Combinatorial Theory \textbf{18} (2011), 2207--2231.

\bibitem{WaltonWhittaker:23}
J.~Walton, M.~Whittaker, \emph{An aperiodic tile with edge-to-edge orientational matching rules},
J. Inst. Math. Jussieu \textbf{22} (2023), no.~4, 1727-1755.

\end{thebibliography}
\end{document}